\documentclass{article}
\usepackage{amssymb,amsmath}
\usepackage{graphicx,CJK}

\def\qed{\hfill {\hbox{${\vcenter{\vbox{               %HOLLOW SQUARE
   \hrule height 0.4pt\hbox{\vrule width 0.4pt height 6pt
   \kern5pt\vrule width 0.4pt}\hrule height 0.4pt}}}$}}}

\def\tr{\triangleright}

\newtheorem{theorem}{Theorem}
\newtheorem{definition}{Definition}

\newtheorem{example}{Example}

\newtheorem{remark}{Remark}

{\qed\par\medskip}

\date{}

\title{\Large \textbf{Kei modules and unoriented link invariants}}

\author{Mike Grier\footnote{Email: mhg02007@mymail.pomona.edu} \and 
Sam Nelson\footnote{Email: knots@esotericka.org}}

\begin{document}
\maketitle

\begin{abstract} We define invariants of unoriented knots and links
by enhancing the integral kei counting invariant $\Phi_X^{\mathbb{Z}}(K)$ for a
finite kei $X$ using representations of the \textit{kei algebra}, 
$\mathbb{Z}_K[X]$, a quotient of the quandle algebra $\mathbb{Z}[X]$ 
defined by Andruskiewitsch and Gra\~{n}a. We give an example that 
demonstrates that the enhanced invariant is stronger than the unenhanced
kei counting invariant. As an application, we use
a quandle module over the Takasaki kei on $\mathbb{Z}_3$ which is not
a $\mathbb{Z}_K[X]$-module to detect the non-invertibility of a virtual 
knot.
\end{abstract}

\medskip

\quad
\parbox{6in}{
\smallskip

}

%\tableofcontents

\section{\large\textbf{Introduction}}

In \cite{T}, Mituhisa Takaski introduced an algebraic structure known
as \textit{kei} (or \begin{CJK*}{UTF8}{min}圭\end{CJK*} in the original kanji).
In \cite{J}
this same structure was reintroduced under the name \textit{involutory quandle},
a special case of a more general algebraic structure related to oriented 
knots known as \textit{quandles}. These algebraic structures can be understood
as arising from the unoriented and oriented Reidemeister moves respectively
via a certain labeling scheme, encoding knot structures in algebra.

In \cite{AG}, for every finite quandle $X$ an associative algebra 
$\mathbb{Z}[X]$ was defined with generators representing coefficients of
``beads'' indexed by quandle labelings of arcs, with relations defined
from the Reidemeister moves. Representations of $\mathbb{Z}[X]$, known 
as \textit{quandle modules}, were used to define new invariants of 
oriented knots and links in \cite{CEGS}. In \cite{HHNYZ} a 
modification of $\mathbb{Z}[X]$ for finite racks (a generalization
of quandles to the case of blackboard-framed isotopy) was used to define 
invariants of framed and unframed oriented knots and links. 

In this paper 
we define a modification of the quandle algebra we call the \textit{kei 
algebra} $\mathbb{Z}_K[X]$ and use it to extend the invariants defined in 
\cite{HHNYZ} to unoriented knots and links. The paper is organized as follows. 
In section \ref{K} we review the basics of kei and the kei counting invariant. 
In section \ref{KM} we define the kei algebra and kei modules. In section
\ref{I} we define the 
kei module enhanced counting invariant. As an application, we use a module
over $\mathbb{Z}[X]$ for a kei $X$ which is not a $\mathbb{Z}_K[X]$-module
to detect the non-invertibility of a virtual knot.
In section \ref{Q} we collect a few questions for future research.

\section{\large\textbf{Kei}}\label{K}

\textit{Kei}  or 
\textit{involutory quandles} were introduced by Mituhisa Takasaki in
1945 \cite{T} and later reintroduced independently by David Joyce 
and S.V. Matveev in the early 1980s \cite{J,M}. 

\begin{definition}
\textup{A \textit{kei} or \textit{involutory quandle} is a set $X$
with a binary operation $\tr$ satisfying for all $x,y,z\in X$
\begin{list}{}{}
\item[(i)] $x\tr x=x$,
\item[(ii)] $(x\tr y)\tr y=x$, and
\item[(iii)] $(x\tr y)\tr z=(x\tr z)\tr(y\tr z)$.
\end{list}}
\end{definition}

\begin{example}
\textup{Let $X$ be any abelian group regarded as a $\mathbb{Z}$-module.
Then $X$ is a kei under the operation}
\[x\tr y = 2y-x.\]
\textup{Such a kei is known as a \textit{Takasaki kei}. If 
$X\cong\mathbb{Z}_n$ then $X$ is often denoted as $R_n$ in the knot theory
literature, known as the \textit{dihedral quandle} on $n$ elements. $R_n$
can also be understood as the set of reflections of a regular $n$-gon.}
\end{example}

\begin{example}
\textup{Let $X$ be any module over $\mathbb{Z}[t]/(t^2-1)$. Then $X$ is a 
kei known as an \textit{Alexander kei} under the operation}
\[x\tr y= tx+(1-t)y.\]
\textup{A Takasaki kei is an Alexander kei with $t=-1.$}
\end{example}

\begin{example}
\textup{Let $L$ be an unoriented link diagram and let $A=\{a_1,\dots,a_n\}$ 
be a set of generators corresponding bijectively with the set of arcs of $L$. 
The \textit{Fundamental Kei} of $L$, $FK(L)$, is defined in the following way. 
First, let $W(L)$ be the set of \textit{kei 
words} in $A$, defined recursively by the rules}
\begin{itemize}
\item $a\in A\Rightarrow a\in W(L)$ \textup{and}
\item $x,y\in W(L)\Rightarrow x\tr y\in W(L)$.
\end{itemize}
\textup{Then the \textit{free kei on A} is the set of equivalence classes of
kei words in $A$ under the equivalence relation generated by relations of the
forms}
\begin{itemize}
\item $x\tr x\sim x$,
\item $(x\tr y)\tr y\sim x$, \textup{and}
\item $(x\tr y)\tr z\sim (x\tr z)\tr (y\tr z)$
\end{itemize}
\textup{for all $x,y,z\in W(L)$. The free kei is a kei under the operation
$[x]\tr [y]=[x\tr y]$. Now, at each crossing in $L$, we have a 
\textit{crossing relation} given by $z=x\tr y$ where $y$ is the overcrossing
arc and $x$ and $z$ are the undercrossing arcs. That is, we have}
\[\includegraphics{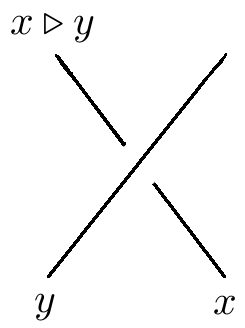}.\]
\textup{Then the \textit{fundamental kei} of $L$, $FK(L)$, is the set of 
equivalence classes of free kei elements modulo the crossing relations
of $L$, or equivalently $FK(L)$ is the set of equivalence classes of 
kei words in $A$ modulo the equivalence relation determined by the
crossing relations together with the free kei relations.} 
\end{example}

It is convenient to describe a finite kei $X=\{x_1,\dots,x_n\}$ with a 
matrix encoding the operation table of $X$, i.e. a matrix $M_X$ whose 
$(i,j)$ entry is $k$ where $x_k=x_i\tr x_j$. For example, the Takasaki
kei on $\mathbb{Z}_3$ has matrix
\[M_X=
\left[\begin{array}{ccc} 
1 & 3 & 2 \\ 
3 & 2 & 1 \\ 
2 & 1 & 3 \\
\end{array}\right]\]
where we set $x_1=0,\ x_2=1$ and $x_3=2$.

As with groups and other algebraic structures, we have the following
standard notions:
\begin{definition} \textup{
Let $X$ and $Y$ be kei.
\begin{itemize}
\item A map $f:X\to Y$ is a \textit{kei homomorphism} if for all
$x,x'\in X$ we have $f(x\tr x')=f(x)\tr f(x')$;
\item A subset $Y\subset X$ which is itself a kei under the kei operation
$\tr$ of $X$ is a \textit{subkei} of $X$. It is easy to check that
$Y\subset X$ is a subkei if and only if $Y$ is closed under $\tr$.
\end{itemize}}
\end{definition}

For defining invariants of unoriented links, we have the following
well-known result:
\begin{theorem}
If $L$ and $L'$ are ambient isotopic unoriented links, then there is
an isomorphism of kei $\phi:FK(L)\to FK(L')$. For any finite kei
$X$, the sets of homomorphisms $\mathrm{Hom}(FK(L),X)$ and
$\mathrm{Hom}(FK(L'),X)$ are finite and 
there is an induced bijection $\phi_{\ast}:\mathrm{Hom}(FK(L),X)\to
\mathrm{Hom}(FK(L'),X)$. In particular, the cardinality 
$\Phi_X^{\mathbb{Z}}(L)=|\mathrm{Hom}(FK(L),X)|$ is a non-negative 
integer-valued invariant of unoriented links known as the 
\textit{integral kei counting invariant}.
\end{theorem}

A kei homomorphism $f:FK(L)\to X$ can be represented as a labeling of
the arcs of $L$ with elements of $X$ satisfying the crossing relations 
at every crossing -- such a labeling defines a unique homomorphism, and
every $f\in \mathrm{Hom}(FK(L),X)$ can be so represented.

\begin{example}
\textup{We can use the kei counting invariant to see that the trefoil
knot $3_1$ is nontrivially knotted. Let $X$ be the Takasaki kei on 
$\mathbb{Z}_3$; we have $x\tr y=2y-x=2y+2x$. The crossing relations
in $3_1$ give us the system of linear equations}
\[\begin{array}{rcl}
z & = & 2x + 2y \\
y & = & 2z + 2x \\
x & = & 2x + 2y
\end{array} \rightarrow
\left[\begin{array}{ccc}
2 & 2 & 2 \\
2 & 2 & 2 \\
2 & 2 & 2 \\
\end{array}\right]\rightarrow
\left[\begin{array}{ccc}
1 & 1 & 1 \\
0 & 0 & 0 \\
0 & 0 & 0 \\
\end{array}\right]
\]
\textup{and the solution space is two-dimensional, giving us a total of 
$\Phi_{X}^{\mathbb{Z}}(3_1)=9$ solutions. Since 
$\Phi_{X}^{\mathbb{Z}}(\mathrm{Unknot})=3$, the integral kei counting invariant
detects the knottedness of the trefoil.}
\end{example}

\begin{remark}
\textup{Replacing the second kei axiom with the alternative axiom
\begin{list}{}{}
\item[(ii')] There exists a second operation $\tr^{-1}$ satisfying
$(x\tr y)\tr^{-1} y=x=(x\tr^{-1} y)\tr y$ for all $x,y\in X$
\end{list}
yields an algebraic object known as a \textit{quandle}, which is the
oriented analog of kei. Labeling oriented
links according to the signed crossing conditions}
\[\includegraphics{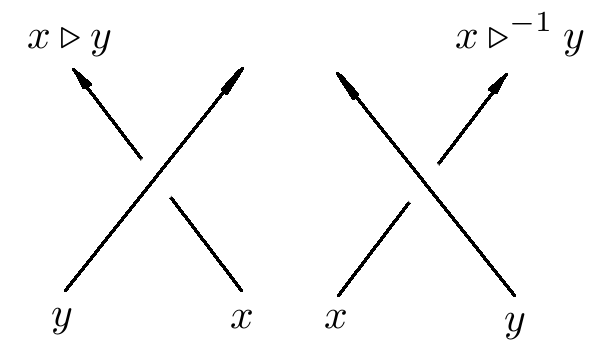}\]
\textup{defines homomorphisms
from the fundamental quandle of the link $L$ into $X$; the \textit{integral
quandle counting invariant} $\Phi_X^{\mathbb{Z}}(L)$ is then an invariant
of oriented links.}
\end{remark}

\section{\large\textbf{Kei algebras and modules}}\label{KM}

Let $X$ be a finite kei. We would like to define an associative algebra
on $X$ generated by ``beads'' such that secondary labelings of $X$-labeled
link diagrams by beads are preserved by Reidemeister moves. Specifically,
at a crossing in a link diagram with arcs labeled $x,y$ and $x\tr y$, we 
define the following relationship between the beads $a,b$ and $c$:
\[
\includegraphics{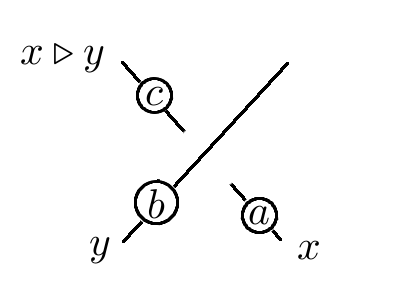} \quad \raisebox{0.5in}{$c=t_{x,y}a+s_{x,y}b.$}
\]

The \textit{kei algebra} of $X$, $\mathbb{Z}_K[X]$, will be the quotient
of the polynomial algebra $\mathbb{Z}[t_{x,y},s_{x,y}]$ by the ideal $I$
required to obtain invariance under unoriented Reidemeister moves. 

First, we note that the bead relationship above also requires that 
$a=t_{x\tr y,y}c+s_{x\tr y,y}b$; together these imply
\[
a=t_{x,y}t_{x\tr y,y}a+(t_{x,y}s_{x\tr y,y}+s_{x,y})b,\]
which yields
\begin{equation}\label{eq1}
t_{x,y}t_{x\tr y,y}=1 \quad \mathrm{and} \quad
t_{x,y}s_{x\tr y,y}+s_{x,y}=0.
\end{equation}

From the Reidemeister I move, we must have $t_{x,x}+s_{x,x}=1$:
\begin{equation}\label{eq2}
\raisebox{-0.9in}{\includegraphics{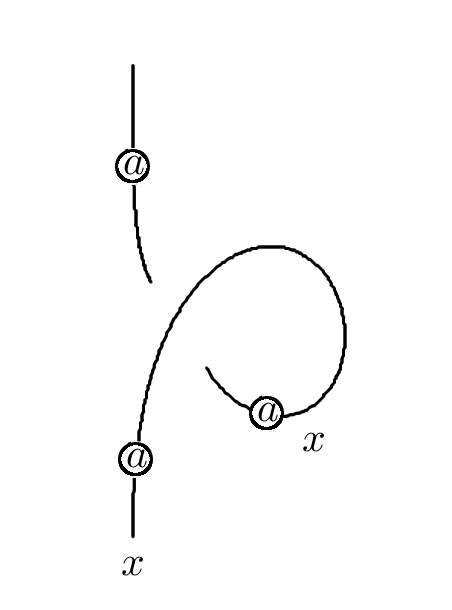}} \quad a=t_{x,x}a+s_{x,x}a
\end{equation}

The Reidemeister II move yields conditions equivalent to equation (\ref{eq1}): 
\[\includegraphics{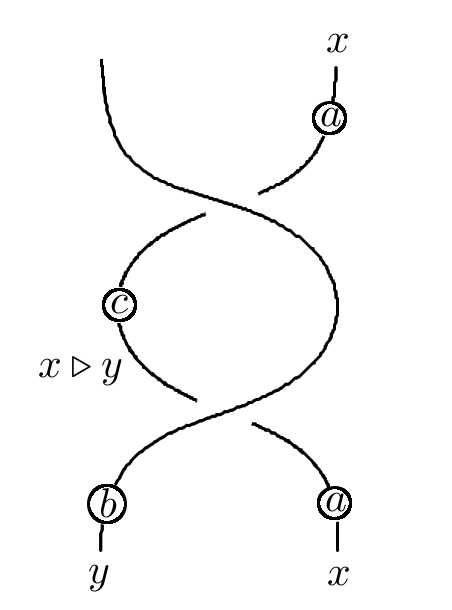} 
\raisebox{1.2in}{$\begin{array}{c}
\begin{array}{rcl}
c & = & t_{x,y}a+s_{x,y}b \\
a & = & t_{x\tr y,y}c+s_{x\tr y,y}b \\
\Rightarrow a & = & t_{x\tr y,y}t_{x,y}a+(t_{x\tr y,y}s_{x,y}+s_{x\tr y,y})b \\
\end{array}\\
\end{array}$}\]

The Reidemeister III move yields the defining equations for the original
rack algebra $\mathbb{Z}[X]$ from \cite{AG}: 
\[\includegraphics{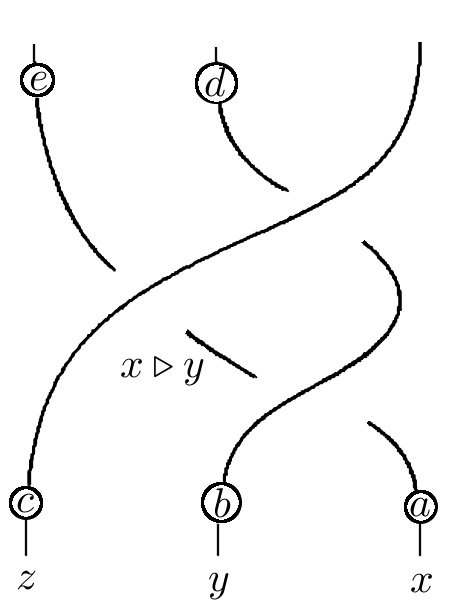} \raisebox{0.9in}{$\sim$}
\includegraphics{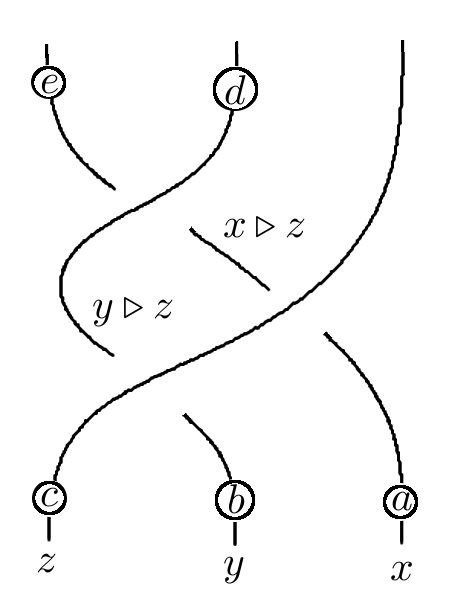} \quad
\raisebox{0.9in}{$\begin{array}{rcl}
e & = & t_{x\tr y,z}t_{x,y}a+t_{x\tr y,z}s_{x.y} b + s_{x\tr y,z} c \\
 & = & t_{x\tr z,y\tr z}t_{x,z}a +s_{x\tr z,y\tr z}t_{y,z}b \\ 
 & & +(t_{x\tr z,y\tr z}s_{x,z}+s_{x\tr z,y\tr z}s_{y,z})c. \\
\end{array}$}
\]
\begin{equation}\label{eq3}
t_{x\tr y,z}t_{x,y}=t_{x\tr z, y\tr z}t_{x,z},
\quad t_{x\tr y,z}s_{x,y}=s_{x\tr z,y\tr z}t_{y,z}\quad
\mathrm{and} \quad s_{x\tr y,z}=s_{x\tr z,y\tr z}s_{y,z}+t_{x\tr z,y\tr z}s_{x,z}.
\end{equation}

We can now define the kei algebra of a finite kei $X$.

\begin{definition}
\textup{Let $X$ be a finite kei. The \textit{kei algebra} $\mathbb{Z}_K[X]$ of 
$X$ is the quotient of the polynomial algebra 
$\mathbb{Z}[t_{x,y},s_{x,y}]$ for all $x,y\in X$ by the ideal $I$ generated by 
all elements of the forms}
\begin{itemize}
\item $t_{x,y}s_{x\tr y,y}+s_{x,y}$,
\item $t_{x,y}t_{x\tr y,y}-1$,
\item $t_{x,x}+s_{x,x}-1$,
\item $t_{x\tr y,z}t_{x,y}-t_{x\tr z, y\tr z}t_{x,z},$
\item $t_{x\tr y,z}s_{x,y}-s_{x\tr z,y\tr z}t_{y,z}$, \textup{and}
\item $s_{x\tr y,z}-s_{x\tr z,y\tr z}s_{y,z}-t_{x\tr z,y\tr z}s_{x,z}$
\end{itemize}
\textup{for all $x,y,z\in X$. A \textit{$\mathbb{Z}_K[X]$-module} or just
an $X$-module is a representation of $\mathbb{Z}_K[X]$, i.e. an abelian
group $A$ with a family of automorphisms $t_{x,y}:A\to A$ and endomorphisms
$s_{x,y}:A\to A$ satisfying the conditions (1), (2) and (3) above.}
\end{definition}

\begin{example}\label{ex1}
\textup{Let $X$ be a kei. Any ring $R$ becomes a $\mathbb{Z}_K[X]$-module
by choosing invertible elements $t_{x,y}$ and elements $s_{x,y}$
for $x,y\in X$ satisfying the conditions (1), (2) and (3). In particular,
if $X=\{x_1,x_2,\dots,x_n\}$ is a finite kei, we can specify a 
$\mathbb{Z}_K[X]$-module structure on $R$ with a $n\times 2n$ block
matrix $M_R=[T|S]$ where $T(i,j)=t_{x_i,x_j}$ and $S(i,j)=s_{x_i,x_j}$.}
\end{example}

\begin{remark}
\textup{The \textit{quandle algebra} defined in \cite{AG} is the quotient
of the polynomial algebra $\mathbb{Z}[t_{x,y}^{\pm 1},s_{x,y}]$ by the ideal 
generated by the relations coming from the Reidemeister I and III moves, 
i.e.,
\begin{itemize}
\item $t_{x,x}+s_{x,x}-1$
\item $t_{x\tr y,z}t_{x,y}-t_{x\tr z, y\tr z}t_{x,z},$
\item $t_{x\tr y,z}s_{x,y}-s_{x\tr z,y\tr z}t_{y,z},$
\item $s_{x\tr y,z}-s_{x\tr z,y\tr z}s_{y,z}-t_{x\tr z,y\tr z}s_{x,z}.$
\end{itemize}
with the type II move condition handled by the bead labeling rule below.}
\[\includegraphics{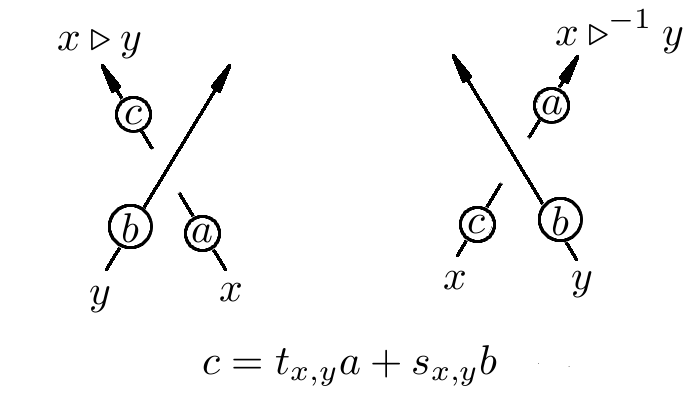}\]
\textup{ The
kei algebra $\mathbb{Z}_K[X]$ is a quotient of the quandle algebra by
the additional relations}
\[ t_{x,y}s_{x\tr y,y}+s_{x,y} \quad \mathrm{and}\quad 1-t_{x,y}t_{x\tr y,y}.
\]
\end{remark}

\begin{example}\label{ex2}
\textup{For a specific instance of the type of kei module defined in 
example \ref{ex1}, let $X$ be the 3-element Takasaki kei with kei matrix}
\[M_X=\left[\begin{array}{ccc}
1 & 3 & 2 \\
3 & 2 & 1 \\
2 & 1 & 3
\end{array}\right]\]
\textup{and let $R=\mathbb{Z}_5$. Our \texttt{python} computations indicate
that there are $48$ $\mathbb{Z}_K[X]$-module structures on $R$, including for
instance}
\[M_R=\left[\begin{array}{ccc|ccc}
4 & 1 & 3 & 2 & 4 & 1 \\
3 & 4 & 2 & 3 & 2 & 3 \\
2 & 1 & 4 & 4 & 1 & 2
\end{array}\right].\]
\end{example}

\begin{remark}\label{rem1}
\textup{For a given kei $X$, the set of $\mathbb{Z}_K[X]$-modules over a 
given ring $R$ is a subset of the set of $\mathbb{Z}[X]$-modules, and can 
be a proper subset depending on $R$, since a $\mathbb{Z}[X]$-module satisfies
the conditions in equation (1) and (3) but not necessarily those of equation
(2). For instance, our \texttt{python} computations reveal a total of $32$ 
$\mathbb{Z}[X]$-modules on the kei $X$ and ring $R$ in example \ref{ex2} 
which are not $\mathbb{Z}_K[X]$-modules, including for instance}
\[M_R=\left[\begin{array}{ccc|ccc}
2 & 1 & 2 & 4 & 2 & 3 \\
1 & 2 & 2 & 2 & 4 & 3 \\
4 & 4 & 2 & 4 & 4 & 4
\end{array}\right].\]
\textup{The invariants defined in the next section associated with such
modules are invariants of oriented links but not invariants of unoriented
links.}
\end{remark}

\begin{example}
\textup{Another important example of a $\mathbb{Z}_K[X]$ module is the
\textit{fundamental $\mathbb{Z}_K[X]$-module of an $X$-labeled link}. 
Let $L$ be an unoriented link with a labeling $f:FK(L)\to X$ by a kei $X$. 
On each arc of $L$, we place a bead; the set of crossing relations then 
determines a presentation for a $\mathbb{Z}_K[X]$-module, denoted 
$\mathbb{Z}_f[X]$,
which we can represent concretely
with a coefficient matrix of the resulting homogeneous system of linear 
equations. For instance, let $X$ be the kei with matrix}
\[M_X=\left[\begin{array}{ccc} 
1 & 1 & 2 \\
2 & 2 & 1 \\
3 & 3 & 3
\end{array}\right];\]
\textup{then the (4,2)-torus link with the $X$-labeling below has fundamental 
$\mathbb{Z}_k[X]$-module presented by the matrix $M_{\mathbb{Z}_f[X]}$:}
\[\includegraphics{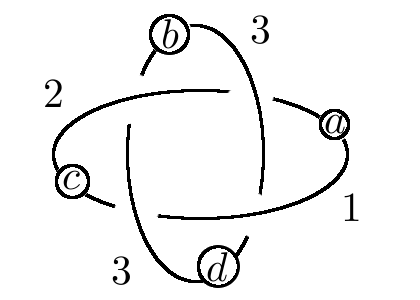}
\quad\raisebox{0.5in}{$M_{\mathbb{Z}_f[X]}=
\left[\begin{array}{rrrr}
t_{13} & s_{13} & -1 & 0 \\
0 & t_{32} & s_{32} & -1 \\
-1 & 0 & t_{23} & s_{23} \\
s_{31} & -1 & 0 & t_{31}
\end{array}\right].$}\]
\end{example}

\section{\large\textbf{Kei module enhancements of the counting invariant}}\label{I}

We can now define invariants of unoriented knots and links using kei modules.
The idea is to use the set of homomorphisms $g:\mathbb{Z}_f[X]\to R$
from the fundamental kei module of an $X$-labeled diagram $L$ to the 
kei module $R$ as a signature for each kei homomorphism $f:FK(L)\to X$.

\begin{definition}
\textup{Let $L$ be an unoriented knot or link, $X$ a finite kei and
$R$ a finite $\mathbb{Z}_K[X]$-module. The \textit{kei module enhanced 
multiset} invariant of $L$ associated to $X$ and $R$ is the multiset of
cardinalities of the sets of $\mathbb{Z}_k[X]$-module homomorphisms, i.e., }
\[\Phi_{X,R}^{K,M}(L)=\left\{
|\mathrm{Hom}_{\mathbb{Z}_K[X]}(\mathbb{Z}_f[X], M)| \ :\ 
f\in \mathrm{Hom}(FK(L),X)
\right\}.\]
\textup{Taking the generating function of this multiset gives us a 
polynomial-form invariant for easy comparison: the \textit{kei module 
enhanced invariant} of $L$ with respect to $X$ and $M$ is}
\[\Phi_{X,R}^{K}(L)=\sum_{f\in \mathrm{Hom}(FK(L),X)} u^{|\mathrm{Hom}_{\mathbb{Z}_K[X]}(\mathbb{Z}_f[X], M)|}.\]
\end{definition}

By construction, we have the following:
\begin{theorem}
If $L$ and $L'$ are ambient isotopic unoriented links, $X$ is a finite kei
and $R$ is a $\mathbb{Z}_K[X]$-module, then $\Phi_{X,R}^{K,M}(L)=
\Phi_{X,R}^{K,M}(L')$ and $\Phi_{X,R}^{K}(L)=\Phi_{X,R}^{K}(L')$.
\end{theorem}

\begin{remark}
\textup{If $R$ is not a finite ring, we can replace the infinite cardinality
$|\mathrm{Hom}_{\mathbb{Z}_K[X]}(\mathbb{Z}_f[X], M)|$ with the rank of the
set $\mathrm{Hom}_{\mathbb{Z}_K[X]}(\mathbb{Z}_f[X], M)$ as a 
$\mathbb{Z}_K[X]$-module.}
\end{remark}

To compute $\Phi_{X,R}^{K}$, for each kei labeling $f:FK(L)\to X$ of $L$ by 
$X$, we first obtain the matrix for $\mathbb{Z}_f[X]$, replace each $t_{x,y}$ 
and $s_{x,y}$ with its value in $R$, and solve the resulting system of equations
to obtain the contributions to $\Phi_{X,R}^{K}$ for $f$.

\begin{example}
\textup{Let $L$ be the figure eight knot $4_1$ and let $X$ and $R$ be the 
kei and kei module on $\mathbb{Z}_5$ from example \ref{ex2}. The set of 
$X$-labelings of $L$ includes only constant labelings, i.e. every arc is
labeled with a 1, 2 or 3. For example, the constant labeling with
every arc labeled 1 yields the listed $\mathbb{Z}_f[X]$-presentation 
matrix:} 
\[\includegraphics{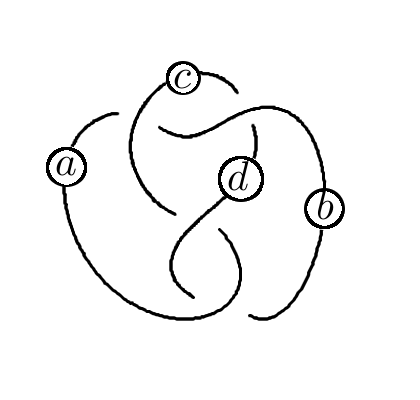}\quad
\raisebox{0.5in}{$M_{\mathbb{Z}_f[X]}=\left[\begin{array}{rrrr}
t_{11} & -1 & s_{11} & 0 \\
0 & s_{11} & t_{11} & -1 \\
-1 & 0 & t_{11} & s_{11} \\
s_{11} & t_{11} & 0 & -1
\end{array}\right]$}\]
\textup{Replacing the $t_{xy}$ and $s_{x,y}$ with their values in $R$ and
row-reducing over $\mathbb{Z}_5$, we obtain}
\[\left[\begin{array}{rrrr}
4 & 4 & 2 & 0 \\
0 & 2 & 4 & 4 \\
4 & 0 & 4 & 2 \\
2 & 4 & 0 & 4
\end{array}\right]
\longrightarrow
\left[\begin{array}{rrrr}
1 & 1 & 3 & 0 \\
0 & 1 & 2 & 2 \\
0 & 0 & 0 & 0 \\
0 & 0 & 0 & 0 \\
\end{array}\right]
\]
\textup{and this $X$-labeling contributes a $u^{25}$ to the invariant
$\Phi_{X,R}^{k}(4_1).$ Summing these contributions over the complete set
of $X$-labelings gives us $\Phi_{X,R}^K(4_1)=3u^{25}$. Comparing this to 
the unknot, which has $\Phi_{X,R}^{k}(\mathrm{Unknot})=3u^5$, we see that
$\Phi_{X,R}^{K}$ distinguishes the unoriented figure eight from the unoriented
unknot despite the two having equal kei counting invariant values. In 
particular, since $\Phi_X^{\mathbb{Z}}(k)$ is obtained from $\Phi_{X,R}^{K}$
by evaluating at $u=1$,  $\Phi_{X,R}^{K}$ 
is a strictly stronger invariant than $\Phi_X^{\mathbb{Z}}(k)$. }
\end{example}

\begin{example}
\textup{Our \texttt{python} computations yield the listed values for
$\Phi_{X,R}^K$ with $X$ the 3-element Takasaki kei and the randomly selected
$\mathbb{Z}_K[X]$-module over $\mathbb{Z}_7$ below for the prime knots with 
up to eight crossings and prime links with up to seven crossings as listed in 
the \texttt{knot atlas} \cite{KA}:}
\[M_R=\left[\begin{array}{rrr|rrr}
6 & 3 & 5 & 2 & 5 & 3 \\
5 & 6 & 3 & 3 & 2 & 5 \\
3 & 5 & 6 & 5 & 3 & 2
\end{array}\right]\]
\[\begin{array}{r|l}
\Phi_{X,M}^{K}(L) & L \\ \hline
3u^7 & \mathrm{unknot}, 4_1, 5_1, 6_2, 6_3, 7_2, 7_3, 7_5, 7_6, 8_1, 8_2, 8_3, 8_4, 8_6, 8_7, 8_8, 8_9, 8_{12} 8_{13} 8_{14}, 8_{17},\\
 &  L2a1, L4a1, L6a2, L6a4, L6n1, L7a2, L7a3, L7a4, L7a7, L7n1, L7n2 \\
3u^7+6u^{49} & 3_1, 6_1, 7_4, 8_{10}, 8_{11}, 8_{15}, 8_{19}, 8_{20}, 8_{21}, L6a1, L6a3, L6a5, L7a1, L7a5 \\
3u^{7}+24u^{49} & 8_{18} \\
3u^{49} & 5_2, 7_1,8_{16}, L7a6 \\ 
9u^{49} & 7_7, 8_5  \\
\end{array}
\]
\end{example}

\begin{remark}
\textup{As with most enhancements of quandle-related counting invariants,
$\Phi_{X,M}^{K}$ is well-defined for unoriented virtual links as well as 
classical links.}
\end{remark}

In our final example, we use a quandle module which is not a 
kei module to detect the non-invertibility of a virtual knot.

\begin{example}
\textup{Let $X$ be the kei from example \ref{ex2} and $M$ the quandle module
from remark \ref{rem1}. Since $M$ is not a kei module, $\Phi_{X}^{M}$ is an 
invariant of oriented knots and links, but not unoriented knots and
links. Thus, we can potentially use $\Phi_X^M$ to compare the two orientations
of a non-invertible knot. In particular, consider the virtual knot numbered
$4.97$ in the Knot Atlas \cite{KA}; it is the closure of the virtual braid
below. Let us denote $4.97$ with the upward orientation by $4.97_{\uparrow}$ and
$4.97$ with the downward orientation as $4.97_{\downarrow}.$
The only labelings of $4.97$ by $X$ are constant labelings, of which 
there are three for both orientations, the unenhanced integral kei 
counting invariant $\Phi_X^{\mathbb{Z}}(4.97_{\uparrow})=3=
\Phi_X^{\mathbb{Z}}(4.97_{\downarrow})$, and $\Phi_X^{\mathbb{Z}}$ does not 
distinguish $4.97_{\uparrow}$
from $4.97_{\downarrow}$. However, the constant labeling with every arc 
labeled with a $1\in X$ yields the listed fundamental kei module presentation 
matrices.
Replacing $t_{1,1}$ and $s_{1,1}$ with their values from $M$ yields the
listed matrices, which we row-reduce over $\mathbb{Z}_5$ to obtain the
cardinalities of the solution spaces which form the signature of the constant 
labeling by the element $1\in X$.}
\[\includegraphics{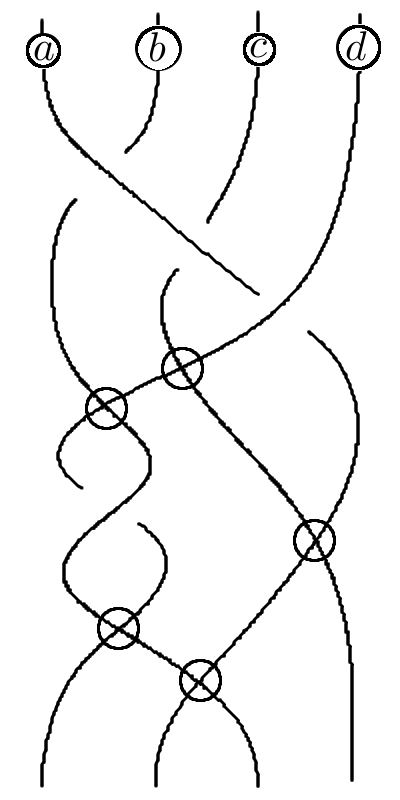}
\raisebox{1.2in}{$\begin{array}{rl}
M_{\mathbb{Z}[f]}(4.97_{\downarrow}): &
\left[\begin{array}{cccc} 
s_{11} & -1 & t_{1,1} & 0 \\
s_{11} & 0 & -1 & t_{1,1} \\
t_{11} & -1 & 0 & s_{1,1} \\
-1 & 0 & s_{11} & t_{1,1} \\
\end{array}\right]
\rightarrow
\left[\begin{array}{cccc} 
4 & 4 & 2 & 0 \\
4 & 0 & 4 & 2 \\
2 & 4 & 0 & 4 \\
4 & 0 & 4 & 2
\end{array}\right]\rightarrow
\left[\begin{array}{cccc} 
1 & 0 & 1 & 3 \\
0 & 1 & 2 & 2 \\
0 & 0 & 0 & 0 \\
0 & 0 & 0 & 0
\end{array}\right] \\
&  \\
&  \\
M_{\mathbb{Z}[f]}(4.97_{\uparrow}): &
\left[\begin{array}{cccc} 
s_{1,1} & t_{1,1} & -1 & 0 \\
s_{1,1} & 0 & t_{1,1} & -1 \\
-1 & t_{1,1} & 0 & s_{1,1} \\
t_{1,1} & 0 & s_{1,1} & -1 \\
\end{array}\right]\rightarrow
\left[\begin{array}{cccc} 
4 & 2 & 4 & 0 \\
4 & 0 & 2 & 4 \\
4 & 2 & 0 & 4 \\
2 & 0 & 4 & 4 \\
\end{array}\right]\rightarrow
\left[\begin{array}{cccc} 
1 & 0 & 0 & 4 \\
0 & 1 & 0 & 4 \\
0 & 0 & 1 & 4 \\
0 & 0 & 0 & 0 \\
\end{array}\right]
\end{array}$}
\]
\textup{Since $t_{1,1}=t_{2,2}=t_{3,3}=2$ and $s_{1,1}=s_{2,2}=s_{3,3}=4$,
we get the same signatures for all three labelings for each knot,
respectively $u^{25}$ and $u^5$, and thus we have}
\[\Phi_X^M(4.97^{\downarrow})=3u^{25}\ne 3u^{5}=\Phi_X^M(4.97_{\uparrow})\]
\textup{and for non-kei module quandle modules $M$ over a finite kei,
$X$, the quandle module enhanced counting invariant $\Phi_X^M$ is capable
of detecting invertibility of virtual (and hence classical) knots.}
\end{example}

\section{\large\textbf{Questions}}\label{Q}

In this section we collect a few open questions for future research.

In our computations we have only considered the simplest type of 
$\mathbb{Z}_K[X]$ modules, namely $\mathbb{Z}_k[X]$-module structures on
$\mathbb{Z}_n$ with the action of $t_{x,y}$ and $s_{x,y}$ given by 
multiplication by fixed elements of $\mathbb{Z}_K[X]$. Expanding to 
other abelian groups and other automorphisms $t_{x,y}:X\to X$ and 
endomorphisms $s_{x,y}:X\to X$ should give interesting results. We
are particularly interested in the case of non-commuting
$t_{x,y}$ and $s_{x,y}$ values.

We have generalized the rack module bead counting invariant from 
\cite{HHNYZ}, but several other oriented link invariants using the 
quandle algebra were defined in \cite{CEGS}; these invariants should 
have generalizations to the unoriented case using the kei algebra.

\bigskip

\noindent
\textsc{Department of Mathematics \\
Pomona College \\
610 N. 6th  \\
Claremont, CA 91711} 

\medskip

\noindent
\textsc{Department of Mathematical Sciences \\
Claremont McKenna College \\
850 Columbia Ave. \\
Claremont, CA 91711}

\end{document}